# A Probability Inequality for Convolutions of MTP$_2$-Distribution Functions

Thomas Royen

TH Bingen, University of Applied Sciences

Berlinstrasse 109, D55411 Bingen, Germany

e-mail: thomas.royen@t-online.de

## Abstract.

A probability inequality is proved for n-fold convolutions of a smooth cumulative distribution function on $\mathbb{R}_+^n$, which is multivariate totally positive of order 2 (MTP$_2$). This inequality is better than an inequality of the same form as the Gaussian correlation inequality for distribution functions. An example are some multivariate chi-square distributions, derived from the diagonal of a Wishart matrix.

## 1. Introduction

If $F \in C^2(\mathbb{R}^n)$ or $F \in C^2(\mathbb{R}_+^n)$ is a strict MTP$_2$ - cumulative distribution function (cdf) then

$$\frac{\partial^2}{\partial x_j \partial x_i} \log F(\mathbf{x}) > 0, \; i \neq j, \qquad (1)$$

holds for all $\mathbf{x} \in \mathbb{R}^n$ or $\mathbf{x} \in \mathbb{R}_+^n = (0,\infty)^n$ respectively. Furthermore, it was proved in [4] for the events $A_i = \{X_i \leq x_i\}$ that

$$\frac{P\left(\bigcap_{i=1}^n A_i\right)}{P\left(\bigcap_{i \in I} A_i\right) P\left(\bigcap_{i \notin I} A_i\right)}, \; \varnothing \neq I \subset \{1,\ldots,n\}, \qquad (2)$$

is decreasing in each $x_i$ (in [4], without the assumption "strict", only "non-increasing"). Consequently, with $B_i = \{X_i \leq x_i\}$, the inequality

$$P\left(\bigcap_{i=1}^n A_i\right) \geq \frac{P\left(\bigcap_{i=1}^n B_i\right)}{P\left(\bigcap_{i \in I} B_i\right) P\left(\bigcap_{i \notin I} B_i\right)} P\left(\bigcap_{i \in I} A_i\right) P\left(\bigcap_{i \notin I} A_i\right) \qquad (3)$$

holds for $\mathbf{x} \leq \mathbf{b}$, i.e. $x_i \leq b_i$, $i=1,\ldots,n$. This inequality is better than an inequality of the same type as the Gaussian correlation inequality (GCI) for cdfs, as long as $P\left(\bigcap_{i=1}^n B_i\right) > P\left(\bigcap_{i \in I} B_i\right) P\left(\bigcap_{i \notin I} B_i\right)$. The general GCI was proved in [6] and in its special form for cdfs extended to some multivariate gamma distributions. Distributions, satisfying the inequalities $P\left(\bigcap_{i=1}^n A_i\right) \geq P\left(\bigcap_{i \in I} A_i\right) P\left(\bigcap_{i \notin I} A_i\right)$, are called too "strongly positive lower orthant dependent" (SPLOD). For some further consequences of the inequality (3) see [4].

In [4] some improved Gaussian correlation inequalities are given, extended to multivariate chi-square and gamma distributions, which were originally derived from the diagonal of a Wishart $W_n(\nu, \Sigma)$-matrix, but these multivariate gamma distributions can be defined by the Laplace transform

---

$$|\mathbf{I}+\mathbf{\Sigma T}|^{-\alpha} \tag{4}$$

with the identity matrix $\mathbf{I}$, the associated covariance matrix $\mathbf{\Sigma}$ (not to confuse with the covariances of the components), $\mathbf{T}=Diag(t_1,...,t_n)$, $t_i \geq 0$, and real "degrees of freedom" $\nu=2\alpha \in \mathbb{N} \cup \left([(n-1)/2],\infty\right)$, see [5]. The $\Gamma_n(\alpha,\mathbf{\Sigma})$- cdf belonging to (4) is here denoted by $G_\alpha(\mathbf{x},\mathbf{\Sigma})$. The most comprehensive collection for integral and series representations and approximations for $\Gamma_n(\alpha,\mathbf{\Sigma})$- cdfs is possibly found in the appendix of [2]. In [4] some examples are found of MTP2 - cdfs $G_\alpha(x,\mathbf{\Sigma})$. The $\Gamma_n(\frac{1}{2},\mathbf{\Sigma})$ - distributed random vector $(X_1,...,X_n)$ is distributed as $\left(\frac{1}{2}Z_1^2,...,\frac{1}{2}Z_n^2\right)$, where $(Z_1,...,Z_n)$ has an $N_n(\mathbf{0},\mathbf{\Sigma})$ - distribution. According to [3], the vector $(|Z_1|,...,|Z_n|)$, and consequently $(X_1,...,X_n)$, is MTP2 iff there exists a sign matrix $\mathbf{S}=Diag(s_1,...,s_n)$, $s_i=\pm 1$, which generates an M-matrix $\mathbf{S\Sigma^{-1}S}$, i.e. with only non-positive off-diagonal elements and with only non-negative off-diagonal elements in $\mathbf{S\Sigma S}$. This is equivalent to the infinite divisibility of the $\Gamma_n(\frac{1}{2},\mathbf{\Sigma})$ - distribution, see [1]. The Laplace transform in (4) shows the invariance of the corresponding $G_\alpha(\mathbf{x},\mathbf{\Sigma})$ - cdf under the transformation $\mathbf{\Sigma} \to \mathbf{S\Sigma S}$. The MTP2 - property is generally not stable under convolutions. The most simple example of an MTP2 - cdf $G_{1/2}(\mathbf{x},\mathbf{\Sigma})$ arises from a non-singular $\mathbf{\Sigma}$ with correlations $\rho_{ij}=a_ia_j$, $i\neq j$, all $a_i \in (-1,1)$. (At most one $|a_i|=1$ is admissible.) For this example the MTP2 - property holds for all $\alpha$ too.

## 2. The Inequality

The proof of the inequality in the following Theorem is based on the decrease of the ratio in (2). To make this paper more self-contained, the short proof for this decrease is given here again.

**Proposition.** If the cdf $F \in C^2(\mathbb{R}^n)$ (or $F \in C^2(\mathbb{R}_+^n)$) is strict MTP2 with an everywhere positive density $f$ on $\mathbb{R}^n$ (or $f$ on $\mathbb{R}_+^n$), then the ratio

$$\frac{P\left(\bigcap_{i=1}^n A_i\right)}{P\left(\bigcap_{i\in I} A_i\right)P\left(\bigcap_{i\notin I} A_i\right)}, \quad \emptyset \neq I \subset \{1,...,n\},$$

is decreasing in every $x_i$.

**Proof.** Let e.g., $i$ be an index of $I$. With the above assumption this decrease is equivalent to

$$\frac{\partial}{\partial x_i}\log P\left(\bigcap_{k=1}^n A_k\right) < \frac{\partial}{\partial x_i}\log P\left(\bigcap_{k\in I} A_k\right). \tag{5}$$

According to formula (1) it is

$$\frac{\partial}{\partial x_j}\left(\frac{\partial}{\partial x_i}\log P\left(\bigcap_{k=1}^n A_k\right)\right) > 0, \ i \in I, j \notin I.$$

Since $\lim_{\min x_j \to \infty, j\notin I} P\left(\bigcap_{k=1}^n A_k\right) = P\left(\bigcap_{k\in I} A_k\right)$, the inequality (5) is satisfied. $\square$

The following probability inequality (6) is a consequence of the decreasing ratio in (2). It is proved inductively by successive convolutions. We write the vector $\mathbf{x}$ as $\mathbf{x}=\mathbf{x}_1 \oplus \mathbf{x}_2$, where $\mathbf{x}_1$ has the components $x_k$, $k \in I_1$, $\emptyset \neq I_1 \subset \{1,...,n\}$. The corresponding marginal cdfs of $F(\mathbf{x})$ are denoted by $F_i(\mathbf{x}_i)$ and the corresponding probability densities by $f_i(\mathbf{x}_i)$, $i=1,2$.

**Theorem.** If the cdf $F \in C^2((0,\infty)^n)$ is strict MTP$_2$ with an everywhere positive density function $f$, then the n-fold convolutions $F^{n*}(\mathbf{x}) = F * ... * F$ satisfy for all $n \in \mathbb{N}$, $n \geq 2$, and all $\mathbf{x} \in (0,\infty)^n$ the inequality

$$F^{n*}(\mathbf{x}) > (r(\mathbf{x}))^n F_1^{n*}(\mathbf{x}_1) F_2^{n*}(\mathbf{x}_2) \text{ with } r(\mathbf{x}) = \frac{F(\mathbf{x})}{F_1(\mathbf{x}_1) F_2(\mathbf{x}_2)} > 1, \tag{6}$$

and consequently the inequality

$$F^{n*}(\mathbf{x}) > (r(\mathbf{b}))^n F_1^{n*}(\mathbf{x}_1) F_2^{n*}(\mathbf{x}_2) \tag{7}$$

too on the region $\mathbf{x} \leq \mathbf{b}$.

**Proof.** It is $F(\mathbf{x}) = r(\mathbf{x}) F_1(\mathbf{x}_1) F_2(\mathbf{x}_2)$ and, because of the decreasing ratio in (2), formula (6) implies

$$F^{(n+1)*}(\mathbf{x}) = \int_{0 < \mathbf{y} \leq \mathbf{x}} F^{n*}(\mathbf{x} - \mathbf{y}) f(\mathbf{y}) dy_1...dy_n > (r(\mathbf{x}))^n \int_{0 \leq \mathbf{z} < \mathbf{x}} F_1^{n*}(\mathbf{z}_1) F_2^{n*}(\mathbf{z}_2) f(\mathbf{x} - \mathbf{z}) dz_1...dz_n =$$

$$(r(\mathbf{x}))^n \int_{0 < \mathbf{z} \leq \mathbf{x}} f_1^{n*}(\mathbf{z}_1) f_2^{n*}(\mathbf{z}_2) F(\mathbf{x} - \mathbf{z}) dz_1...dz_n >$$

$$(r(\mathbf{x}))^{n+1} \int_{0 < \mathbf{z} \leq \mathbf{x}} f_1^{n*}(\mathbf{z}_1) f_2^{n*}(\mathbf{z}_2) F_1(\mathbf{x}_1 - \mathbf{z}_1) F_2(\mathbf{x}_2 - \mathbf{z}_2) dz_1...dz_n =$$

$$(r(\mathbf{x}))^{n+1} F_1^{(n+1)*}(\mathbf{x}_1) F_2^{(n+1)*}(\mathbf{x}_2). \quad \square$$

A similar inequality can be derived for convolutions of independent not identically distributed MTP$_2$-random vectors on $(0,\infty)^n$, but then the powers $(r(\mathbf{x}))^n$ have to be replaced by products $\prod_{i=1}^n r_i(\mathbf{x})$.

For MTP$_2$- gamma-cdfs $G_{1/2}(\mathbf{x},\mathbf{\Sigma})$ we obtain directly from the above theorem with

$$\mathbf{\Sigma} = \begin{pmatrix} \mathbf{\Sigma}_{11} & \mathbf{\Sigma}_{12} \\ \mathbf{\Sigma}_{21} & \mathbf{\Sigma}_{22} \end{pmatrix}, \ rank(\mathbf{\Sigma}_{12}) > 0, \text{ and } r(\mathbf{x}) = r(\tfrac{1}{2},\mathbf{x},\mathbf{\Sigma}) = \frac{G_{1/2}(\mathbf{x},\mathbf{\Sigma})}{G_{1/2}(\mathbf{x}_1,\mathbf{\Sigma}_{11}) G_{1/2}(\mathbf{x}_2,\mathbf{\Sigma}_{22})} > 1$$

the inequality

$$G_\alpha(\mathbf{x},\mathbf{\Sigma}) > (r(\tfrac{1}{2},\mathbf{x},\mathbf{\Sigma}))^{2\alpha} G_\alpha(\mathbf{x}_1,\mathbf{\Sigma}_{11}) G_\alpha(\mathbf{x}_2,\mathbf{\Sigma}_{22}) \text{ with } 2\alpha \in \mathbb{N}, 2\alpha \geq 2, \tag{8}$$

and consequently

$$G_\alpha(\mathbf{x},\mathbf{\Sigma}) > (r(\tfrac{1}{2},\mathbf{b},\mathbf{\Sigma}))^{2\alpha} G_\alpha(\mathbf{x}_1,\mathbf{\Sigma}_{11}) G_\alpha(\mathbf{x}_2,\mathbf{\Sigma}_{22}) \text{ for all } \mathbf{x} \leq \mathbf{b}. \tag{9}$$

For larger degrees of freedom $\nu = 2\alpha$ this inequality is coarse, but nevertheless remarkably better than the inequality for only SPLOD. The ratio $r(\tfrac{1}{2},\mathbf{x},\mathbf{\Sigma})$ is (apart from simulations) numerically available for rather high dimensions by means of suitable software for multivariate Gaussian distributions. In a similar way it follows

$$G_{\alpha+k/2}(\mathbf{x},\mathbf{\Sigma}) > (r(\tfrac{1}{2},\mathbf{x},\mathbf{\Sigma}))^k G_{\alpha+k/2}(\mathbf{x}_1,\mathbf{\Sigma}_{11}) G_{\alpha+k/2}(\mathbf{x}_2,\mathbf{\Sigma}_{22}) \tag{10}$$

with an MTP$_2$– cdf $G_{1/2}(\mathbf{x},\mathbf{\Sigma})$, $k \in \mathbb{N}$ and any real values $2(\alpha+1) > [(n-1)/2]$, see page 5 in [4].